\documentclass{article}
\usepackage{graphicx}
\usepackage{latexsym}
\usepackage{amssymb}
\usepackage{amsmath}
\usepackage{layout}
\newtheorem{teo}{Theorem}[section]
\newtheorem{lema}[teo]{Lemma}

\newenvironment{dem}{\text\bf {\it Proof:}}{}
\newcommand{\R}{\mathbb{R}}
\newcommand{\N}{\mathbb{N}}
\newcommand{\eps}{\varepsilon}
\begin{document}

\begin{center}
{\huge {\ On the convergence to the multiple
Wiener-It\^{o} integral} }

\vspace{.5cm}

\normalsize {\bf Xavier Bardina$^{1,*}$, Maria Jolis$^{1}$}
 and {\bf Ciprian A. Tudor$^2$}

\vspace{.25cm}

{\footnotesize \it $^1$Departament de Matem\`atiques; Universitat
Aut\`onoma de Barcelona; 08193-Bellaterra (Barcelona); Spain.;
bardina@mat.uab.cat; mjolis@mat.uab.cat

 $^2$SAMOS/MATISSE, Centre d'Economie de La Sorbonne, Universit\'e de Panth\'eon-Sorbonne Paris 1, 90, rue de Tolbiac,
75634 Paris Cedex 13, France.; tudor@univ-paris1.fr

*Corresponding author}

\end{center}

%
\vspace{0.5cm}

\begin{abstract}
We study the convergence to the multiple Wiener-It\^{o} integral
from processes with absolutely continuous paths. More precisely,
consider a family of processes, with paths in the Cameron-Martin
space, that converges weakly to a standard Brownian motion in
$\mathcal C_0([0,T])$. Using these processes, we construct a family
that converges weakly, in the sense of the finite dimensional
distributions, to the multiple Wiener-It\^{o} integral process of a
function $f\in L^2([0,T]^n)$. We prove also the weak convergence in
the space $\mathcal C_0([0,T])$ to the second order integral for two
important families of processes that converge to a standard Brownian
motion.

\end{abstract}

\vskip0.5cm

{\bf  2000 AMS Classification Numbers: } 60B10, 60F05,  60H05.

 \vskip0.3cm


\vskip0.3cm

\section{Introduction and preliminaries}

Let $Y$ be a semimartingale with trajectories belonging to the space
$\mathcal D([0,T])$ of functions right continuous with left limits
in all point in $[0,T]$, and define the following iterated It\^o
integrals
\[J_k(Y)_t=\left\{\begin{array}{ll}
                  Y_t  &  \mbox{if $k=1$}\\
                  \int_0^t J_{k-1}(Y)_{s^{-}} \, dY_s &\mbox{for $k\geq2.$}
                  \end{array}
            \right. \]

\noindent Suppose that $\{X^{\eps}\}_{\eps>0}$ is a family of
semimartingales wih paths in $\mathcal D([0,T])$ that converges
weakly in this space to another semimartingale $X$, as $\eps$ tends
to zero. Avram (1988) proved that the following statements are
equivalent

\begin{eqnarray*}
\mathcal
L(X^\varepsilon,[X^\varepsilon,X^\varepsilon])&\stackrel{w}\longrightarrow
&\mathcal L(X,[X,X])\qquad {\rm when }
\quad\varepsilon\downarrow0,\qquad {\rm and} \cr && \cr \mathcal
L(J_1(X^\varepsilon),\dots,J_m(X^\varepsilon))&\stackrel{w}\longrightarrow
&\mathcal L(J_1(X),\dots,J_m(X))\qquad {\rm when
}\quad\varepsilon\downarrow0,
\end{eqnarray*}
(
here $\stackrel{w}\longrightarrow$ denotes the weak convergence in $\mathcal D([0,T])^{2} $ and $\mathcal D([0,T])^{m}$
respectively) where, if we denote by $Y^c$ the continuous part of a semimartingale $Y$, the process $[Y,Y]$ is defined by
$$[Y,Y]_{\,t}=<Y^c,Y^c>_t+\sum_{s\le t}(\Delta Y_s)^2.$$

This result shows that in order to obtain (joint) weak convergence
of It\^{o} multiple integrals we need the convergence of $X^{\eps}$ to
$X$ but also also the convergence of the second order variations.
When our semimartingale is the Wiener process, there is a lot of
important examples of families of processes with absolutely
continuous paths converging in law in $\mathcal C([0,T])$ to it. In
this case, clearly, we do not have convergence of the quadratic
variations to that of the Brownian motion.

 Consider the Cameron-Martin space:
$$\mathcal H:=\{\eta\in\mathcal C([0,T]):\,\eta_t=\int_{0}^t \eta_s' ds,\,\eta'\in L^2([0,T])\},$$
 and a family of
processes $(\eta_{\varepsilon})_{\varepsilon>0}$ with paths
belonging to the Cameron Martin space given by
\begin{equation}\label{etas}
\eta_{\varepsilon}(t)=\int_0^t \theta_{\varepsilon}(s) \,ds,
\end{equation}
 such that $(\eta_{\eps})_{\eps>0}$ converges weakly to a standard
 Brownian motion in $\mathcal C_0([0,T])$, the space of continuous
 functions defined in $[0,T]$ which are null at zero.

 Consider now, for a function $f\in L^2([0,T]^n)$, the multiple
 integrals
 $$I_{\eta_{\varepsilon}}(f)_t=\int_0^t\cdots\int_0^t
f(t_1,\dots,t_n)\,d\eta_{\varepsilon}(t_1)\cdots
d\eta_{\varepsilon}(t_n).$$ In Bardina and Jolis (2000) the
convergence in law of $(I_{\eta_{\varepsilon}}(f))_{\eps}$ was
studied. The authors proved that in order to obtain convergence for
all families  $(\eta_{\varepsilon})$ with values in $\mathcal H$ and
converging in law to the Wiener process, the function $f$ needs to
be given by a multimeasure. For other classes of functions, only
partial results were obtained with some particular families
$(\eta_{\varepsilon})$. In all the cases the limit was the
Stratonovich integral of $f$ with respect to the Wiener process, as
defined by Sol\'{e} and Utzet (1990). This fact is not surprising,
taking into account that the multiple Stratonovich integral must
satisfy the rules of the ordinary calculus. On the other hand, this
integral is a complicated object, it is defined by a limiting
procedure and only some classes of functions (as tensor products or
continuous functions) are recognized as Stratonovich integrable.

A natural question is that of the possibility of obtaining, for  a function $f$, its multiple Wiener-It\^{o}-type integral as
a limit in law of some multiple integrals  with respect to the absolutely continuous processes $\eta_{\eps}$. Since in the
definition of the multiple Wiener-It\^{o} integral with respect to the Wiener process, the approximating procedure implies the
suppression of  the values on the diagonals, one can expect that a similar idea will allow to obtain this integral as a limit
law.

We denote  by $Y^f_{\eta_{\varepsilon}}$ the stochastic processes
defined by
\begin{align}\label{familia}
Y_{\eta_{\varepsilon}}^f(t):=\int_{[0,t]^n}f(x_1,x_2,\dots,x_n)\prod^n_{\footnotesize{\begin{array}{c}i,j=1\\i\neq
j\end{array}}}
I_{\{|x_i-x_j|>\varepsilon\}}d\eta_{\varepsilon}(x_1)\cdots\cdots
d\eta_{\varepsilon}( x_n)\nonumber\\
=\int_{[0,t]^n}f(x_1,x_2,\dots,x_n)\prod_{i=1}^n\theta_{\varepsilon}(x_i)\prod^n_{\footnotesize{\begin{array}{c}i,j=1\\i\neq
j\end{array}}}I_{\{|x_i-x_j|>\varepsilon\}}dx_1\cdots
dx_n,\end{align} for all $t\in [0,T]$. We have studied the weak
 convergence of the finite dimensional distributions of $(Y^f_{\eta_{\varepsilon}})_{\eps>0}$ to that of the
corresponding multiple Wiener-It\^{o}-integral of $f$ with respect
to the Wiener process, and also the convergence in $\mathcal
C_0([0,T])$  of second order integrals for two important families of
process $(\eta_{\eps})$. With regard to the convergence of the
finite dimensional distributions, we have proved that there is
convergence under rather general conditions on $(\eta_{\eps})$, see
Theorem 2.3. For the convergence in  $\mathcal C_0([0,T])$ of the
second order integral, we have proved it for the so-called Donsker
and Kac-Stroock families of processes. It is worth to note that in
all the results obtained here, the function $f$ is an arbitrary
function in $L^2([0,T]^n)$, that is, all the domain of the
Wiener-It\^{o} integral. This is a very different situation from
that of Bardina and Jolis (2000).

Section 2 deals with the problem of the convergence of finite dimensional distributions and Section 3 is devoted to prove
convergence in the space of continuous functions for the second order integral with respect to the Donsker and Kac-Stroock
processes. In all the paper we denote the positive multiplicative constants that do not depend neither on $\eps$ nor on the
function $f$ by $C$, although their values can change from an expression to another one.

\section{Convergence of the finite dimensional distributions}
\subsection{Some general results}
We first state a general lemma that will be the main tool in order
to prove the convergence of the finite dimensional distributions. We
state it in our particular setting.
\begin{lema}\label{lema1} Let $(S,\|\cdot\|)$ be a normed space and consider
$\{J^{\varepsilon}\}_{\varepsilon\geq0}$ a family of linear
applications defined on $S$ with values in the space of
m-dimensional finite a.s. random variables,
$\left(L^0(\Omega)\right)^m$.

\noindent Denote by $|\cdot|$ the Euclidian norm in $\R^m$.

\noindent Assume that there exists a positive constant $C$ such that
for all $f\in S$
\begin{equation}\label{H}
\sup_{\varepsilon\geq0}E\left|J^{\varepsilon}(f)\right|\leq C\|f\|.
\end{equation}

\noindent Assume also that there exists a dense subset $D\subset S$
such that for all $f\in D$, $J^{\varepsilon}(f)$ converges in law to
$J^0(f),$ when $\varepsilon$ tends to 0.

\noindent Then, $J^{\varepsilon}(f)$ converges in law to $J^0(f),$
for all $f\in S$, when $\varepsilon$ tends to 0.
\end{lema}

We will denote by ${\cal{E}}^{\prime ,n}$ the space of simple functions on $[0,T]^n$ that can be written as
\begin{equation}\label{simpesp}
f(x_1,\dots,x_n)=\sum_{k=1}^m\alpha_k
I_{\Delta_k}(x_1,\dots,x_n),
 \end{equation}
 where, $m\in\mathbb N$, $\alpha_k\in\R$, for all
$k\in\{1,\dots,m\}$, and
$\Delta_k=(a_k^1,b_k^1]\times(a_k^2,b_k^2]\times\cdots\times(a_k^n,b_k^n]$
with $[a_k^h,b_k^h]\cap[a_k^l,b_k^l]=\emptyset$ for all $h\neq l$.

\begin{lema}\label{lema2}
Let $\left(\eta_{\varepsilon}\right)_{\varepsilon>0}$ be a family of
processes with trajectories in the Cameron-Martin space $\mathcal H$
of the form (\ref{etas}). Assume that the finite dimensional
distributions of the family
$\left(\eta_{\varepsilon}\right)_{\varepsilon>0}$ converge in law to
the finite dimensional distributions of a standard Brownian motion
$W$ when $\varepsilon$ tends to 0.

\noindent Consider $f\in{\cal{E}}^{\prime ,n}$. Then, the finite dimensional distributions of the processes
$Y_{\eta_{\varepsilon}}^f$ defined in (\ref{familia}) converge in law to the finite dimensional distributions of the multiple
Wiener-It\^{o} integral $I_n(f\cdot I_{[0,t]^n})$ when $\varepsilon$ tends to 0.
\end{lema}

\begin{dem} Since $f\in{\cal{E}}^{\prime ,n}$
it follows that, for $\varepsilon$ small enough,
$$f(x_1,x_2,\dots,x_n)\prod^n_{\footnotesize{\begin{array}{c}i,j=1\\i\neq j\end{array}}}
I_{\{|x_i-x_j|>\varepsilon\}}=f(x_1,x_2,\dots,x_n),$$
for all $(x_1,x_2,\dots,x_n)\in[0,T]^n$. And then, if $f$ is given
by (\ref{simpesp}),  for those $\varepsilon$,
$$Y_{\eta_{\varepsilon}}^f(t)=\sum_{k=1}^m\alpha_k\prod_{i=1}^n\left[\eta_{\varepsilon}(b_k^i\wedge t)-\eta_{\varepsilon}(a_k^i\wedge t)\right].$$

We conclude from
the convergence of the finite dimensional distributions of $\eta
_{\varepsilon}$ to that of the Brownian motion that for all
$t_1,\dots,t_r\in[0,T]$ the vector
$(Y_{\eta_{\varepsilon}}^f(t_1),\dots,Y_{\eta_{\varepsilon}}^f(t_r))$
converges in law to
$$
\left( \sum_{k=1}^m\alpha_k\prod_{i=1}^n\left[W(b_k^i\wedge
t_1)-W(a_k^i\wedge t_1)\right],\dots,
\sum_{k=1}^m\alpha_k\prod_{i=1}^n\left[W(b_k^i\wedge
t_r)-W(a_k^i\wedge t_r)\right]
 \right),
$$
when $\varepsilon$ tends to 0.

But, since $f\in{\cal{E}}^{\prime ,n}$, by the definition of the
multiple Wiener-It\^{o} integral (see It\^{o}, 1951) the last random
vector equals to
$$\left(I_n(f\cdot I_{[0,t_1]^n}),\dots,I_n(f\cdot I_{[0,t_r]^n})\right).$$

\hfill $\Box$

\end{dem}

\bigskip

The following theorem is the main result of this section and gives sufficient conditions for the family $(\eta_{\varepsilon})$
in order to have the convergence of the finite dimensional distributions of $Y^f_{\eta_{\varepsilon}}$ to those of the
multiple Wiener-It\^{o} integral process, for any $f\in L^2([0,T]^n)$.
\begin{teo}\label{prop1}
Let $\left(\eta_{\varepsilon}\right)_{\varepsilon>0}$ be a family of processes with trajectories in the Cameron-Martin space
$\mathcal H$ of the form (\ref{etas}). Assume that the finite dimensional distributions of the family
$\left(\eta_{\varepsilon}\right)_{\varepsilon>0}$ converge in law to the finite dimensional distributions of a standard
Brownian motion when $\varepsilon$ tends to 0.

\noindent Assume also that there exists a positive constant $C$ such
that \begin{equation} \label{ccc}\sup_{\varepsilon>0,t\in[0,T]}
E\left|Y_{\eta_{\varepsilon}}^f(t)\right|\leq C\|f\|_{L^2([0,T]^n)},
\end{equation}for all $f\in L^2([0,T]^n)$.

\noindent Then, the finite dimensional distributions of the family of processes $\{Y_{\eta_{\varepsilon}}^f\}_{\varepsilon>0}$
converge in law to those of the multiple Wiener-It\^{o} integral $I_n(f\cdot I_{[0,t]^n})$ for all $f\in L^2([0,T]^n)$, when
$\varepsilon$ tends to 0.
\end{teo}

\begin{dem}
Take $t_1,\dots,t_r\in[0,T]$. In order to see that for all $f\in
L^2([0,T]^n)$, the random vector
$(Y_{\eta_{\varepsilon}}^f(t_1),\dots,Y_{\eta_{\varepsilon}}^f(t_r))$
converges in law to
$$\left(I_n(f\cdot I_{[0,t_1]^n}),\dots,I_n(f\cdot I_{[0,t_r]^n})\right),$$
when $\varepsilon$ tends to 0, we will apply Lemma \ref{lema1}. Take
$S= L^2([0,T]^n)$  and consider, for all $\varepsilon>0$,  the
linear operators
\begin{eqnarray*}
J^{\varepsilon}:L^2([0,T]^n)&\longrightarrow&\left(L^0(\Omega)\right)^r\\
f&\longrightarrow&(Y_{\eta_{\varepsilon}}^f(t_1),\dots,Y_{\eta_{\varepsilon}}^f(t_r)),
\end{eqnarray*}
and, for $\varepsilon=0$, the linear operator,
\begin{eqnarray*}
J^{0}:L^2([0,T]^n)&\longrightarrow&\left(L^0(\Omega)\right)^r\\
f&\longrightarrow&\left(I_n(f\cdot I_{[0,t_1]^n}),\dots,I_n(f\cdot
I_{[0,t_r]^n})\right).
\end{eqnarray*}

Condition (\ref{H}) of Lemma \ref{lema1} is satisfied because, by
hypothesis,
$$\sup_{\varepsilon>0} E\left|J^{\varepsilon}(f)\right|\leq C\|f\|_{L^2([0,T]^n)},$$
and, on the other hand, it is well known that
$$E\left|J^{0}(f)\right|\leq C\|f\|_{L^2([0,T]^n)}.$$

By Lemma \ref{lema2} we have  that, for all $f\in{\cal{E}}^{\prime ,n}$, $J^{\varepsilon}(f)$ converges in law to $J^{0}(f)$.
This fact completes the proof because ${\cal{E}}^{\prime ,n}$ is a dense subset of $L^2([0,T]^n)$.

\hfill $\Box$

\end{dem}

\bigskip


We can also consider the  problem of the vectorial convergence to
multiple Wiener-It\^o integrals  in the sense of the finite
dimensional distributions. Fix a natural number $d\geq 2 $ and
consider $d$ integers $n_{1}, n_{2}, \ldots , n_{d}\geq 1$. Let
$f_{k}\in L^{2} \left( [0,T]^{n_{k}}\right) $ for $ k=1, \ldots d$
and consider the sequence of stochastic processes with values in
$\mathbb{R}^{d}$
\begin{equation}
\label{Zeps} Z^{\varepsilon }(t) =\left(
Y_{\eta_{\varepsilon}}^{f_{1} }(t), \ldots ,
Y_{\eta_{\varepsilon}}^{f_{d}} (t)\right), \hskip0.5cm t\in [0,T].
\end{equation}
with $Y_{\eta_{\varepsilon}}^{f_{k}}$, $k=1, \ldots d,$ defined by
(\ref{familia}). We can prove the next theorem that shows the
convergence, as $\varepsilon \to 0$, of the finite dimensional
distributions of $Z^{\varepsilon}$ to those of the vector of
multiple Wiener-It\^o integrals
\begin{equation}
\label{Z} Z(t)=  \left(  I_{n_{1}}(f_{1} \cdot I_{[0,t]^{n_{1}}}), \ldots ,  I_{n_{d}}(f_{d} \cdot I_{[0,t]^{n_{d}}})\right),
\hskip0.5cm t\in [0,T].
\end{equation}

\begin{teo}
Let $\left( \eta _{\varepsilon} \right) _{\varepsilon >0} $ be a
family of stochastic processes with paths in ${\cal{H}}$ that
converges in the sense of the finite dimensional distributions to a
standard Brownian motion. Let, for $k=1,\ldots d$,  $f_{k}\in  L^{2}
\left( [0,T]^{n_{k}}\right) $ and assume that condition (\ref{ccc})
is satisfied by every $n_{k}$, $k=1,\ldots ,d$. Then the finite
dimensional distributions of the vector $Z^{\varepsilon}$ given by
(\ref{Zeps}) converges as $\varepsilon \to 0$ to those of the vector
$Z$ given by (\ref{Z}).

\end{teo}
\begin{dem}
The proof follows similar arguments to that of Theorem \ref{prop1}
and then omitted.

\hfill $\Box$
\end{dem}

\subsection{Examples}\label{exa}
We will give now two examples of families $\eta_{\varepsilon}$ for whose  the above Theorems  can be applied.
\subsubsection{Convergence for the Donsker kernels}

Consider now the particular case
$$\theta_{\varepsilon}(s):=\frac{1}{\varepsilon}\sum_{k=1}^{\infty}
 \xi_k I_{[k-1,k)}\left(\frac{s}{\varepsilon^2}\right),$$
 where $\{\xi_k\}$ is a sequence of independent, identically distributed
 random variables satisfying $E(\xi_1)=0$ and $\mbox{Var}(\xi_1)=1$.

 The processes  $\theta_{\varepsilon}$ will be called Donsker
 kernels, because the convergence in law of
  $\eta_{\varepsilon}(t)=\int_0^t\theta_{\varepsilon}(s)ds$ to the
  Brownian motion in $\mathcal C([0,T])$ is given by the well-known Donsker's Invariance Principle.

In view of Theorem \ref{prop1}, in order to prove the convergence of
the finite dimensional distributions of
$\{Y_{\eta_{\varepsilon}}^f\}_{\varepsilon}$ to the finite
dimensional distributions of $I_n(f\cdot I_{[0,t]^n})$, it is enough
to prove that there exists some constant $C>0$ such that, for all
$f\in L^2([0,T]^n)$
$$\sup_{\varepsilon>0,t\in[0,T]}E\Big[\int_{[0,t]^n}f(x_1,x_2,\dots,x_n)
\prod_{i=1}^n\theta_{\varepsilon}(x_i)\!\!\prod^n_{\footnotesize{\begin{array}{c}i,j=1\\i\neq
 j\end{array}}}I_{\{|x_i-x_j|>\varepsilon\}}dx_1\cdots dx_n\Big]^2\leq C\|f\|^2_{L^2([0,T]^n)}.$$

We can assume, without loss of generality that $f$ is symmetric.
Notice that,
\begin{eqnarray}
&&E\Big[\int_{[0,t]^n}f(x_1,x_2,\dots,x_n)\prod_{i=1}^n
\theta_{\varepsilon}(x_i)\prod^n_{\footnotesize{\begin{array}{c}i,j=1\nonumber\\i\neq
j\end{array}}}I_{\{|x_i-x_j|>\varepsilon\}}dx_1\cdots
dx_n\Big]^2\label{estrella}\nonumber\\
&=&\int_{[0,t]^{2n}}f(x_1,x_2,\dots,x_n)f(y_1,y_2,\dots,y_n)E[\prod_{i=1}^n\theta_{\varepsilon}
(x_i)\theta_{\varepsilon}(y_i)]\nonumber\\&\times&
\prod^n_{\footnotesize{\begin{array}{c}i,j=1\\i\neq
j\end{array}}}I_{\{|x_i-x_j|>\varepsilon\}}I_{\{|y_i-y_j|>\varepsilon\}}dx_1\cdots
dx_n dy_1\cdots dy_n.
\end{eqnarray}

We can also suppose that $\varepsilon<1$. In this case, the
condition $|x-y|>\varepsilon$ implies that  $|x-y|>\varepsilon^2$
and then,
\begin{eqnarray*}
&&E\Big[\prod_{i=1}^n\theta_{\varepsilon}(x_i)\theta_{\varepsilon}(y_i)\Big]\prod^n_{\footnotesize{\begin{array}
{c}i,j=1\\i\neq
j\end{array}}}I_{\{|x_i-x_j|>\varepsilon\}}I_{\{|y_i-y_j|>\varepsilon\}}\\
&=&E\Big[\frac1{\varepsilon^{2n}}\!\!\!\!\!\!\sum_{\scriptsize{\begin{array}{c}i_1,\dots,i_n,j_1,\dots,j_n\\i_k\neq
i_l, j_k\neq j_l \, \forall k\neq
l\end{array}}}\!\!\!\!\!\!\prod_{k=1}^n\!\xi_{i_k}\xi_{j_k}I_{[i_k-1,i_k]}
\left(\frac{x_k}{\varepsilon^2}\right)I_{[j_k-1,j_k]}
\left(\frac{y_k}{\varepsilon^2}\right)\Big]
\prod^n_{\scriptsize{\begin{array}{c}i,j=1\\i\neq
j\end{array}}}\!\!\!I_{\{|x_i-x_j|>\varepsilon\}}I_{\{|y_i-y_j|>\varepsilon\}}.
\end{eqnarray*}
Notice that the number of different indexes in each summand
appearing in the above expression is greater or equal than $n$.
Therefore, using also the symmetry of $f$, we can write
(\ref{estrella}) as
\begin{eqnarray*}
&&\int_{[0,t]^{2n}}f(x_1,x_2,\dots,x_n)f(y_1,y_2,\dots,y_n)E\Big[\frac1{\varepsilon^{2n}}{\textstyle\sum'}
\prod_{k=1}^n\xi_{i_k}\xi_{j_k}I_{[i_k-1,i_k]}\left(\frac{x_k}{\varepsilon^2}\right)I_{[j_k-1,j_k]}
\left(\frac{y_k}{\varepsilon^2}\right)\Big]\\
&\times&\prod^n_{\footnotesize{\begin{array}{c}i,j=1\\i\neq
j\end{array}}}I_{\{|x_i-x_j|>\varepsilon\}}I_{\{|y_i-y_j|>\varepsilon\}}dx_1\cdots
dx_1\cdots dx_n dy_1\cdots dy_n\\
&+&\int_{[0,t]^{2n}}f(x_1,x_2,\dots,x_n)f(y_1,y_2,\dots,y_n)E\Big[\frac{n!}{\varepsilon^{2n}}
\!\!\!\!\!\!\sum_{\footnotesize{\begin{array}{c}i_1,\dots,i_n\\i_k\neq
i_l,\, \forall k\neq l\end{array}}}\!\!\!\!\!\!
\prod_{k=1}^n\xi^2_{i_k}I_{[i_k-1,i_k]^2}\left(\frac{x_k}{\varepsilon^2},\frac{y_k}{\varepsilon^2}\right)\Big]\\
&\times&\prod^n_{\footnotesize{\begin{array}{c}i,j=1\\i\neq
j\end{array}}}I_{\{|x_i-x_j|>\varepsilon\}}I_{\{|y_i-y_j|>\varepsilon\}}dx_1\cdots
dx_1\cdots dx_n dy_1\cdots dy_n,
\end{eqnarray*}
where $\sum'$ denotes the sum over all the indexes satisfying that
at least $n+1$ among the $i_1,\dots,i_n$, $j_1,\dots,j_n$ are
different.

 Using now that $\{\xi_k\}$ is a sequence of independent,
identically distributed
 random variables with $E(\xi_1)=0$ and $\mbox{Var}(\xi_1)=1$, we
can bound the last expression by
\begin{eqnarray*}
&&n!\int_{[0,t]^{2n}}|f(x_1,x_2,\dots,x_n)f(y_1,y_2,\dots,y_n)|\\&&\times
\frac{1}{\varepsilon^{2n}}
\left(\!\!\!\!\sum_{\footnotesize{\begin{array}{c}i_1,\dots,i_n\\i_k\neq
i_l,\, \forall k\neq l\end{array}}}\!\!\!\!\!\!
\prod_{k=1}^nI_{[i_k-1,i_k]^2}\left(\frac{x_k}{\varepsilon^2},\frac{y_k}{\varepsilon^2}\right)\right)
dx_1\cdots dx_n dy_1\cdots dy_n\\
&\leq&
n!\int_{[0,t]^{n}}f^2(x_1,x_2,\dots,x_n)\\&&\times\frac{1}{\varepsilon^{2n}}\left(
\!\!\!\!\sum_{\footnotesize{\begin{array}{c}i_1,\dots,i_n\\i_k\neq
i_l,\, \forall k\neq l\end{array}}}\!\!\!\!\!\!
\prod_{k=1}^nI_{[i_k-1,i_k]}\left(\frac{x_k}{\varepsilon^2}\right)\int_{[0,t]^n}\prod_{k=1}^nI_{[i_k-1,i_k]}\left(\frac{y_k}{\varepsilon^2}\right)dy_1\cdots
dy_n\right)
dx_1\cdots dx_n\\
&\le&n!\int_{[0,t]^{n}}f^2(x_1,x_2,\dots,x_n)\left(
\!\!\!\!\sum_{\footnotesize{\begin{array}{c}i_1,\dots,i_n\\i_k\neq
i_l,\, \forall k\neq l\end{array}}}\!\!\!\!\!\!
\prod_{k=1}^nI_{[i_k-1,i_k]}\left(\frac{x_k}{\varepsilon^2}\right)\right)
dx_1\cdots dx_n\\
&\leq&n!\|f\|^2_{L^2([0,T]^n)}.
\end{eqnarray*}

\subsubsection{Convergence for the Kac-Stroock kernels}\label{KSK}

Consider now the following kernels introduced by Kac (1956)
$$\theta_{\varepsilon}(x):=\frac1{\varepsilon}(-1)^{N\left(\frac{x}{\varepsilon^2}\right)},$$
where $N=\{N(s);\,s\geq0\}$ is a standard Poisson process. Stroock
(1982) proved that the family $(\eta_{\varepsilon})_{\varepsilon>0}$
with $\eta_{\varepsilon}(t)=\int_0^t\theta_{\varepsilon}(s)ds$
converges in law in $\mathcal C_0([0,T])$ to the Brownian motion.

As for the Donsker kernels, to prove the convergence of the finite
dimensional distributions of
$\{Y_{\eta_{\varepsilon}}^f\}_{\varepsilon}$ to those of $I_n(f\cdot
I_{[0,t]^n})$, it is enough to prove that there exists some constant
$C>0$ such that, for all $f\in L^2([0,T]^n)$
$$\sup_{\varepsilon>0,t\in[0,T]}E\Big[\int_{[0,t]^n}f(x_1,x_2,\dots,x_n)\prod_{i=1}^n\theta_{\varepsilon}(x_i)
\prod^n_{\footnotesize{\begin{array}{c}i,j=1\\i\neq
j\end{array}}}I_{\{|x_i-x_j|> \varepsilon\}}dx_1\cdots
dx_n\Big]^2\leq C\|f\|^2_{L^2([0,T]^n)}.$$

 Observe
that, denoting by $\mathcal P_n$ the group of permutations of the
set $\{1,\ldots,n\}$, we have that
\begin{eqnarray}\label{laprimera}
&&E[\prod_{i=1}^n\theta_{\varepsilon}(x_i)\theta_{\varepsilon}(y_i)]\prod^n_{\footnotesize{\begin{array}{c}i,j=1\\i\neq
j\end{array}}}I_{\{|x_i-x_j|>\varepsilon\}}I_{\{|y_i-y_j|>\varepsilon\}}\nonumber\\
&=&E[\prod_{i=1}^n\theta_{\varepsilon}(x_i)\theta_{\varepsilon}(y_i)]\!\!\!\prod^n_{\footnotesize{\begin{array}{c}i,j=1\\i\neq
j\end{array}}}\!\!\!I_{\{|x_i-x_j|>\varepsilon\}}I_{\{|y_i-y_j|>\varepsilon\}}\sum_{\sigma,\psi\in\mathcal
P_n}I_{\{x_{\sigma_1}\leq x_{\sigma_2}\leq\cdots\leq
x_{\sigma_n}\}}I_{\{y_{\psi_1}\leq y_{\psi_2}\leq\cdots\leq
y_{\psi_n}\}}\nonumber\\
&=&E[\prod_{i=1}^n\theta_{\varepsilon}(x_i)\theta_{\varepsilon}(y_i)]\!\!\!\prod^n_{\footnotesize{\begin{array}{c}i,j=1\\i\neq
j\end{array}}}\!\!\!I_{\{|x_i-x_j|>\varepsilon\}}I_{\{|y_i-y_j|>\varepsilon\}}\sum_{\sigma,\psi\in\mathcal
P_n} I_{\{x_{\sigma_1},y_{\psi_1}\}\leq \{x_{\sigma_2},y_{\psi_2}\}\leq \cdots\leq \{x_{\sigma_n},y_{\psi_n}\}}\nonumber\\
&&+E[\prod_{i=1}^n\theta_{\varepsilon}(x_i)\theta_{\varepsilon}(y_i)]\!\!\!\prod^n_{\footnotesize{\begin{array}{c}i,j=1\\i\neq
j\end{array}}}\!\!\!I_{\{|x_i-x_j|>\varepsilon\}}I_{\{|y_i-y_j|>\varepsilon\}}\sum_{\sigma,\psi\in\mathcal
P_n}A(x_1,\ldots,x_n,y_1,\ldots,y_n;\sigma,\psi)
\end{eqnarray}
where $\{a,b\}\leq\{c,d\}$ means that $a\vee b\leq c\wedge d$,  and
where $A(x_1,\ldots,x_n,y_1,\ldots,y_n;\sigma,\psi)$ is the sum of
the indicator functions with all the other possible orders between
the $2n$ variables $\{x_{\sigma_1}\leq x_{\sigma_2}\leq\cdots\leq
x_{\sigma_n}\}$ and $\{y_{\psi_1}\leq y_{\psi_2}\leq\cdots\leq
y_{\psi_n}\}$.

We will start with the first summand of the righthand side of
(\ref{laprimera}). Notice that

\begin{eqnarray*}
&E[\prod_{i=1}^n\theta_{\varepsilon}(x_i)\theta_{\varepsilon}(y_i)]I_{\{x_{\sigma_1},y_{\psi_1}\}\leq
\{x_{\sigma_2},y_{\psi_2}\}\leq \cdots\leq
\{x_{\sigma_n},y_{\psi_n}\}}\phantom{xxxxxxxxxxxxxxxxxxx} \\
&=\frac1{\varepsilon^{2n}}E\left[(-1)^{\sum_{i=1}^n\left(N\left(\frac{x_{\sigma_i}}{\varepsilon^2}\right)
+N\left(\frac{y_{\psi_i}}{\varepsilon^2}\right)
\right)}\right]I_{\left\{\{x_{\sigma_1},y_{\psi_1}\}\leq
\{x_{\sigma_2},y_{\psi_2}\}\leq \cdots\leq
\{x_{\sigma_n},y_{\psi_n}\}\right\}}.
\end{eqnarray*}

Using that for $a,\,b\in\N\cup \{0\}$ we have
$(-1)^{a+b}=(-1)^{a-b}$, the fact that the Poisson process has
independent increments, and that if $Z\sim \text{Poiss}(\lambda)$
then $E[(-1)^Z]=\exp(-2\lambda)$, we obtain that the expectation
appearing in the last expression is equal to
$$\exp\left(-2\sum_{i=1}^n\left(\frac{|x_{\sigma_i}-y_{\psi_i}|}{\varepsilon^2}\right)\right).$$

\noindent Moreover,
\begin{eqnarray*}
&&\int_{[0,t]^{2n}}|f(x_1,x_2,\dots,x_n)f(y_1,y_2,\dots,y_n)|\frac1{\varepsilon^{2n}}\exp\left(-2\sum_{i=1}^n\left(\frac{|x_{\sigma_i}-y_{\psi_i}|}{\varepsilon^2}\right)\right)
dx_1\cdots dx_ndy_1\cdots dy_n\\
&\leq&\int_{[0,t]^{n}}f^2(x_1,x_2,\dots,x_n)\frac1{\varepsilon^{2n}}
\left(\int_{[0,t]^n}\exp\left(-2\sum_{i=1}^n
\left(\frac{|x_{\sigma_i}-y_{\psi_i}|}{\varepsilon^2}\right)\right)dy_1\cdots dy_n\right)dx_1\cdots dx_n\\
&\leq&\int_{[0,t]^{n}}f^2(x_1,x_2,\dots,x_n)dx_1\cdots dx_n\\
&=&\|f\|^2_{L^2([0,T]^n)}.
\end{eqnarray*}

We consider now the second summand of (\ref{laprimera}). We have
showed that in the computation of the expectation is important the
order of the $2n$ variables $\{x_{\sigma_1}\leq
x_{\sigma_2}\leq\cdots\leq x_{\sigma_n}\}$ and $\{y_{\psi_1}\leq
y_{\psi_2}\leq\cdots\leq y_{\psi_n}\}$. If we take the variables
$(x_1,\ldots,x_n,y_1\ldots,y_n)$ in each summand of
$A(x_1,\ldots,x_n,y_1\ldots,y_n;\sigma,\psi)$ in groups of two
variables taking into account their order, necessarily one of the
groups will be formed by two variables $x_k,\, x_l$ (for some $k\neq
l\in\{1,\dots,n\}$). Then, when we compute the expectation the
corresponding term will be $$\exp\left(
-2\frac{|x_k-x_l|}{{\varepsilon^2}}\right)$$ and we have that
$$\frac{1}{\varepsilon^{2n}}\exp\left(
-2\frac{|x_k-x_l|}{{\varepsilon^2}}\right)I_{\{|x_k-x_l|>\varepsilon\}}\leq
\frac{1}{\varepsilon^{2n}}e^{-\frac2{\varepsilon}}\leq C.$$

So, we have that
\begin{eqnarray*}
&&\sum_{\sigma,\psi\in\mathcal
P_n}\int_{[0,t]^{2n}}|f(x_1,x_2,\dots,x_n)f(y_1,y_2,\dots,y_n)
|E[\prod_{i=1}^n\theta_{\varepsilon}(x_i)\theta_{\varepsilon}(y_i)]\\&&\times
A(x_1,\ldots,x_n,y_1\ldots,y_n;\sigma,\psi)
\prod^n_{\footnotesize{\begin{array}{c}i,j=1\\i\neq
j\end{array}}}\!\!\!I_{\{|x_i-x_j|>\varepsilon\}}I_{\{|y_i-y_j|>\varepsilon\}}dx_1\cdots
dx_n dy_1\cdots dy_n\\
&\leq&C\int_{[0,t]^{2n}}|f(x_1,x_2,\dots,x_n)f(y_1,y_2,\dots,y_n)|dx_1\cdots
dx_n dy_1\cdots dy_n\\
&\leq&C\|f\|^2_{L^2([0,T]^n)}.
\end{eqnarray*}

\section{Convergence in law in $\mathcal C_0([0,T])$ of the second order integral}

In this section we will see that in the case of the second order
Wiener-It\^o integral, for the examples introduced in Subsection
\ref{exa}, we can prove also the convergence in law in $\mathcal
C_0([0,T]$.

Let us first mention that clearly for every $\varepsilon >0$ the paths of the process $Y_{\eta _{\varepsilon}}^{f}$ are
absolute continuous functions. On the other hand, since the multiple Wiener-It\^o integrals can be expressed as iterate
integrals
$$I_{2}(f\cdot I_{[0,t]^2})=2 \int_{0}^{t} \int_{0}^{y}f(x,y)dW(x) dW(y) $$
for any $f\in L^{2}([0,T]^{2})$, the stochastic process $\left(
I_{2}(f\cdot I_{[0,t]^2})\right) _{t\geq 0}$ admits a version with
continuous trajectories.

When $n=2$, the processes $Y_{\eta_{\varepsilon}}^f$ become
\begin{equation}\label{cas2}Y_{\eta_{\varepsilon}}^f(t):=\int_0^t\int_0^t f(x,y)
\theta_{\varepsilon}(x)\theta_{\varepsilon}(y)I_{\{|x-y|>\varepsilon\}}dxdy,\end{equation}
where $\theta_{\varepsilon}$ are the Kac-Stroock or the Donsker
kernels. In this section we need more integrability for the
variables $\{\xi_k\}$ appearing in the Donsker kernels. Concretely
we will assume that
 $E(\xi_k)^4<+\infty$.

\begin{teo}\label{nouteo} Let $f\in L^2([0,T]^2)$. Then, the processes $Y_{\eta_{\varepsilon}}^f$ given by
(\ref{cas2}) converge weakly to the multiple Wiener-It\^{o} integral of order 2, $I_2(f\cdot I_{[0,t]^2})$, in the space
$\mathcal C_0([0,T])$ when $\varepsilon$ tends to zero.
\end{teo}

\begin{dem}
We can assume without loss of generality that $f$ is symmetric.
  We have proved, in the previous section, the convergence of the finite dimensional
 distributions. So, to prove the convergence in law, it is enough to
 prove that the family of laws of $\{Y_{\eta_{\varepsilon}}^f\}_{\varepsilon}$ is tight in $\mathcal C_0([0,T])$.

It suffices to show that for $s\leq t$
\begin{equation}\label{Billingsley}
E\left(Y_{\eta_{\varepsilon}}^f(t)-Y_{\eta_{\varepsilon}}^f(s)\right)^4\leq
C\left(\int_{[0,T]^2}\bar f^2(x,y)dxdy\right)^2,\end{equation} where
$$\bar f(x,y):=f(x,y)I_{[0,t]^2}(x,y)-f(x,y)I_{[0,s]^2}(x,y).$$

Indeed, for $s\le t$
$$\left(I_{[0,t]^2}-I_{[0,s]^2}\right)^2=I_{[0,t]^2}-I_{[0,s]^2},$$
Therefore, if (\ref{Billingsley}) is satisfied
\begin{eqnarray*}
E\left(Y_{\eta_{\varepsilon}}^f(t)-Y_{\eta_{\varepsilon}}^f(s)\right)^4&\leq&
C\left(\int_{[0,T]^2}\bar f^2(x,y)dxdy\right)^2\\
&=&C\left(\int_{[0,T]^2}f^2(x,y)\left(I_{[0,t]^2}-I_{[0,s]^2}\right)dxdy\right)^2\\
&=&C\left(\int_s^t\int_0^yf^2(x,y)dxdy+\int_s^t\int_0^xf^2(x,y)dydx\right)^2\\
&=&C\left(\int_s^t\int_0^yf^2(x,y)dxdy\right)^2,
\end{eqnarray*}
using the symmetry of $f$ in the last step. Then, by Billingsley
criterium (see Theorem 12.3 of Billingsley (1968)), we will obtain
tightness.

Notice that
\begin{eqnarray}
&&E\left(Y_{\eta_{\varepsilon}}^f(t)-Y_{\eta_{\varepsilon}}^f(s)\right)^4\nonumber\\&=&\int_{[0,T]^8}\prod_{i=0}^3\bar
f(u_{2i+1},u_{2i+2})I_{\{|u_{2i+1}-u_{2i+2}|>\varepsilon\}}
E\left(\prod_{i=1}^8\theta_\varepsilon(u_i)\right)du_1\cdots du_8\nonumber\\
&\leq&\int_{[0,T]^8}\prod_{i=0}^3|\bar
f(u_{2i+1},u_{2i+2})|I_{\{|u_{2i+1}-u_{2i+2}|>\varepsilon\}}\left|E\left(\prod_{i=1}^8\theta_\varepsilon(u_i)\right)\right|du_1\cdots
du_8. \label{eq1}
\end{eqnarray}

$ $From now on  we will study separately the Kac-Stroock case and
the Donsker case.

\subsection*{Proof of Theorem \ref{nouteo} for the Kac-Stroock kernels}
In order to simplify notation denote by $f^S$ the function defined
as
$$f^S(u_1,\dots,u_8)=\sum_{\sigma\in\mathcal P_8}\prod_{i=0}^3|\bar
f(u_{\sigma_{2i+1}},u_{\sigma_{2i+2}})|I_{\{|u_{\sigma_{2i+1}}-u_{\sigma_{2i+2}}|>\varepsilon\}}.$$

In the case of the Kac-Stroock kernels, using the same kind of
computations that in Subsection \ref{KSK}, and using also the
symmetry of $f^{S}$, we have that (\ref{eq1}) can be bounded by
$$C\int_{[0,T]^8}\frac1{\varepsilon^8}f^S(u_1,\dots,u_8)I_{\{u_1<u_2<\dots<u_8\}}\prod_{i=0}^3\exp\left(\frac{-2(u_{2i+2}-u_{2i+1})}{\varepsilon^2}\right).$$

Consider now the different summands appearing in the definition of
$f^S$. Notice that if in a summand appears a factor of the type
$$\exp\left(\frac{-2(x-y)}{\varepsilon^2}\right)I_{\{x-y>\varepsilon\}}$$
we have that
$$\frac1{\varepsilon^8}\exp\left(\frac{-2(x-y)}{\varepsilon^2}\right)I_{\{x-y>\varepsilon\}}\leq\frac1{\varepsilon^8}e^{-\frac2{\varepsilon}}\leq C.$$

And so, all the terms with this type of factors can be bounded by
\begin{eqnarray*}
C\int_{[0,T]^8}\prod_{i=0}^3|\bar f(u_{2i+1},u_{2i+2})|du_1\cdots
du_8&=& C\left(\int_{[0,T]^2}|\bar f(x,y)|dxdy\right)^4\\
&\leq&C\left(\int_{[0,T]^2}\bar f^2(x,y)dxdy\right)^2.
\end{eqnarray*}

For the rest of summands appearing in $f^S$, we bound all the
indicators by 1, and excepting symmetries, there are only two
possible situations:

\subsubsection*{Situation 1} We have terms of the type
\begin{eqnarray*}
&&\int_{[0,T]^8}\frac1{\varepsilon^8}\prod_{i=1}^4|\bar
f(x_{i},y_{i})|\exp\left(\frac{-2|x_1-x_2|}{\varepsilon^2}+\frac{-2|y_1-y_2|}{\varepsilon^2}\right)
\\&&\times\exp\left(\frac{-2|x_3-x_4|}{\varepsilon^2}+\frac{-2|y_3-y_4|}{\varepsilon^2}\right)
 dx_1\dots dx_4 dy_1\dots dy_4.
\end{eqnarray*}
This kind of terms can be bounded  by
\begin{eqnarray*}
&& \int_{[0,T]^8}\frac1{\varepsilon^8}\bar f^2(x_1,y_1)\bar
f^2(x_3,y_3)\exp\left(\frac{-2|x_1-x_2|}{\varepsilon^2}+\frac{-2|y_1-y_2|}{\varepsilon^2}\right)
\\&&\times\exp\left(\frac{-2|x_3-x_4|}{\varepsilon^2}+\frac{-2|y_3-y_4|}{\varepsilon^2}\right)
 dx_1\dots dx_4 dy_1\dots dy_4\\
+&&\int_{[0,T]^8}\frac1{\varepsilon^8}\bar f^2(x_2,y_2)\bar
f^2(x_4,y_4)\exp\left(\frac{-2|x_1-x_2|}{\varepsilon^2}+\frac{-2|y_1-y_2|}{\varepsilon^2}\right)
\\&&\times\exp\left(\frac{-2|x_3-x_4|}{\varepsilon^2}+\frac{-2|y_3-y_4|}{\varepsilon^2}\right)
 dx_1\dots dx_4 dy_1\dots dy_4.
\end{eqnarray*}

Integrating, in the first summand of the last expression, with
respect to $x_2,y_2,x_4,y_4$ and in the second one with respect to
$x_1,y_1,x_3,y_3$ we have that the last expression is bounded by
\begin{eqnarray*}
&&C\int_{[0,T]^4}\bar f^2(x_1,y_1)\bar f^2(x_3,y_3) dx_1dy_1dx_3dy_3
+
C\int_{[0,T]^4}\bar f^2(x_2,y_2)\bar f^2(x_4,y_4) dx_2dy_2dx_4dy_4\\
&=& C\left(\int_{[0,T]^2}\bar f^2(x,y) dxdy\right)^2.
\end{eqnarray*}

\subsubsection*{Situation 2} We have also terms of the type
\begin{eqnarray*}
&&\int_{[0,T]^8}\frac1{\varepsilon^8}\prod_{i=1}^4|\bar
f(x_{i},y_{i})|\exp\left(\frac{-2|x_1-x_2|}{\varepsilon^2}+\frac{-2|x_3-x_4|}{\varepsilon^2}\right.
\\&&+\left.\frac{-2|y_1-y_3|}{\varepsilon^2}+\frac{-2|y_2-y_4|}{\varepsilon^2}\right)
 dx_1\dots dx_4 dy_1\dots dy_4.
\end{eqnarray*}
All these terms are bounded by
\begin{eqnarray*}
&&C \int_{[0,T]^8}\frac1{\varepsilon^8}\bar f^2(x_1,y_1)\bar
f^2(x_4,y_4)\exp\left(\frac{-2|x_1-x_2|}{\varepsilon^2}+\frac{-2|x_3-x_4|}{\varepsilon^2}\right.
\\&&+\left.\frac{-2|y_1-y_3|}{\varepsilon^2}+\frac{-2|y_2-y_4|}{\varepsilon^2}\right)
 dx_1\dots dx_4 dy_1\dots dy_4\\
+&&C\int_{[0,T]^8}\frac1{\varepsilon^8}\bar f^2(x_2,y_2)\bar
f^2(x_3,y_3)\exp\left(\frac{-2|x_1-x_2|}{\varepsilon^2}+\frac{-2|x_3-x_4|}{\varepsilon^2}\right.
\\&&+\left.\frac{-2|y_1-y_3|}{\varepsilon^2}+\frac{-2|y_2-y_4|}{\varepsilon^2}\right)
 dx_1\dots dx_4 dy_1\dots dy_4.
\end{eqnarray*}

Integrating now, in the first summand of the last expression, with
respect to $x_2,y_2,x_3,y_3$ and in the second one with respect to
$x_1,y_1,x_4,y_4$ we have that the last expression is bounded by
\begin{eqnarray*}
&&C\int_{[0,T]^4}\bar f^2(x_1,y_1)\bar f^2(x_3,y_3) dx_1dy_1dx_3dy_3
+
C\int_{[0,T]^4}\bar f^2(x_2,y_2)\bar f^2(x_4,y_4) dx_2dy_2dx_4dy_4\\
&=& C\left(\int_{[0,T]^2}\bar f^2(x,y) dxdy\right)^2.
\end{eqnarray*}


\subsection*{Proof of Theorem \ref{nouteo} for the Donsker kernels}

Remember that in this case
$$\theta_{\varepsilon}(s):=\frac{1}{\varepsilon}\sum_{k=1}^{\infty}
 \xi_k I_{[k-1,k)}(\frac{s}{\varepsilon^2}),$$
 where $\{\xi_k\}$ is a sequence of independent, identically distributed
 random variables satisfying $E(\xi_1)=0$, $\mbox{Var}(\xi_1)=1$ and
 $E(\xi_1)^4<+\infty$.

 Remember also that we can assume that
 $\varepsilon<1$, and then the condition $|x-y|>\varepsilon$ implies that
$|x-y|>\varepsilon^2$.

Expression (\ref{eq1})  equals to
\begin{eqnarray}\label{sumes}
&&\int_{[0,T]^8}\prod_{i=1}^4|\bar
f(u_i,v_i)|I_{\{|u_i-v_i|>\varepsilon^2\}}|E(\prod_{i=1}^4\theta_{\varepsilon}(u_i)\theta_{\varepsilon}(v_i))|du_1\dots
dv_4\nonumber\\
&\phantom{xxxxx}=&\int_{\tiny{[0,T]^8}}\prod_{i=1}^4|\bar
f(u_i,v_i)|I_{\{|u_i-v_i|>\varepsilon^2\}}\nonumber\\
&\phantom{xxxxx}\times& \Big|E\Big(\!\!\!\!\sum_{\tiny{\begin{array}{c}i_1,\dots,i_4 \\
j_1,\ldots,
j_4\end{array}}}\!\!\!\!\xi_{i_1}\!\cdots\!\xi_{j_4}\,I_{\scriptsize{[i_1-1,
i_1]}}(\tiny{\frac{u_1}{\varepsilon^2}})\cdots
I_{\scriptsize{[j_4-1,
j_4]}}(\tiny{\frac{v_4}{\varepsilon^2}})\Big)\Big|\,d u_1\cdots
dv_4.
\end{eqnarray}

We have, on one hand, that the random variables $\xi_i$ are
independent with $E(\xi_k)=0$ and, on the other hand that
$$I_{\{|u_i-v_i|>\varepsilon^2\}}I_{[k-1,k)}\left(
\frac{u_i}{\varepsilon^2}\right)I_{[k-1,k)}\left(\frac{v_i}{\varepsilon^2}\right)=0.$$

Consequently, to compute the expectation in  expression
(\ref{sumes}), we have to consider the different decompositions of
$8$ as sums of natural numbers between 2 and 4: (2+2+2+2), (2+2+4),
(3+3+2) and (4+4), that will be the exponents of the $\xi_i$ in the
products appearing in (\ref{sumes}) with no null expectation. Taking
into account that
$$\sum_{k}I_{\scriptsize{[k-1,
k]}}(\tiny{\frac{u}{\varepsilon^2}})I_{\scriptsize{[k-1,
k]}}(\tiny{\frac{v}{\varepsilon^2}})\le
I_{\{\,|u-v|\,<\,\varepsilon^2\,\}},
$$
that $E(\xi_i^4)<\infty$, and doing similar computations to those of
the last section, the expressions obtained with all these
decompositions, except with the third one (3+3+2), can be bounded by

$$\frac{C}{\varepsilon^8}\int_{[0,T]^8}\prod_{i=0}^3|\bar
f(u_{{2i+1}},u_{{2i+2}})|I_{\{|u_{{2i+1}}-u_{{2i+2}}|>\varepsilon^2\}}
\sum_{\sigma\in\mathcal
P_8}\prod_{i=0}^3I_{\{|u_{\sigma_{2i+1}}-u_{\sigma_{2i+2}}|<\varepsilon^2\}}
du_1\dots du_8.
$$

Observe that the products
$$\prod_{i=0}^3|\bar
f(u_{{2i+1}},u_{{2i+2}})|I_{\{|u_{{2i+1}}-u_{{2i+2}}|>\varepsilon^2\}}
\prod_{i=0}^3I_{\{|u_{\sigma_{2i+1}}-u_{\sigma_{2i+2}}|<\varepsilon^2\}}$$
equal to zero for all permutation $\sigma\in\mathcal P_8$ such that
at least one of the sets of two variables
$\,\{u_1,u_2\},\{u_3,u_4\},\{u_5,u_6\},\{u_7,u_8\}\,$ is transformed
by $\sigma$ in one of them. Then, one can only consider the
permutations $\sigma$ for which, given
$\,\{u_1,u_2\},\{u_3,u_4\},\{u_5,u_6\},\{u_7,u_8\}\,$, there exists
always two couples among them  such that their four variables  are
not paired in the product
$$\prod_{i=0}^3I_{\{|u_{\sigma_{2i+1}}-u_{\sigma_{2i+2}}|<\varepsilon^2\}}.$$
Now, we can proceed as with the Kac-Stroock kernels. If, for
instance, the two couples with the above property are $\{u_1,u_2\}$
and $\{u_3,u_4\}$ we majorize the product
$$\prod_{i=0}^3|\bar
f(u_{{2i+1}},u_{{2i+2}})|\,I_{\{|u_{{2i+1}}-u_{{2i+2}}|>\varepsilon^2\}}$$
by
$$\frac12\Big(\bar f^2(u_1,u_2)\bar f^2(u_3,u_4)+\bar f^2(u_5,u_6)\bar
f^2(u_7,u_8)\big)$$ and, for each summand, perform the integral
first with respect to the remaining four variables. This allows to
cancell the term, $\frac1{\varepsilon^8}$ and we obtain the desired
bound.

Finally, we must to study the term corresponding to the
decomposition (3+3+2). Taking now into account that
$$\sum_{k}I_{\scriptsize{[k-1,
k]}}(\tiny{\frac{u}{\varepsilon^2}})I_{\scriptsize{[k-1,
k]}}(\tiny{\frac{v}{\varepsilon^2}})I_{\scriptsize{[k-1,
k]}}(\tiny{\frac{w}{\varepsilon^2}})\le
I_{\{\,\text{GD}\{u,v,w\,\}\,<\,\varepsilon^2\,\}},
$$
where we denote by $\text{GD}$ the greatest distance between a
sequence of elements, the product of indicators that we will obtain
in this case can be bounded by
$$\sum_{\sigma\in\mathcal P_8}I_{\{\text{GD}\{u_{\sigma_1},u_{\sigma_2},u_{\sigma_3}\}<\varepsilon^2\}}
I_{\{\text{GD}\{u_{\sigma_4},u_{\sigma_5},u_{\sigma_6}\}<\varepsilon^2\}}
I_{\{|u_{\sigma_7}-u_{\sigma_8}|<\varepsilon^2\}}.
$$

Therefore, all the terms, excepting symmetries, will be of the form
\begin{eqnarray*}
&&f(x,y)f(s,t)f(u,v)f(z,w)I_{\{|x-y|>\varepsilon^2\}}I_{\{|s-t|>\varepsilon^2\}}I_{\{|u-v|>\varepsilon^2\}}I_{\{|z-w|>\varepsilon^2\}}\\
&\times&I_{\{\text{GD}\{x,s,u\}<\varepsilon^2\}}
I_{\{\text{GD}\{y,t,z\}<\varepsilon^2\}}
I_{\{|v-w|<\varepsilon^2\}}.
\end{eqnarray*}

(Observe that we do not consider
$I_{\{\text{GD}\{x,s,u\}<\varepsilon^2\}}
I_{\{\text{GD}\{y,t,v\}<\varepsilon^2\}}$ because in this case we
obtain a factor
$I_{\{|z-w|<\varepsilon^2\}}I_{\{|z-w|>\varepsilon^2\}}\,=\,0\,$).

This kind of term can be bounded by
$$f^2(x,y)f^2(z,w)I_A+f^2(s,t)f^2(u.v)I_B,$$
where
$$A:={\{|v-w|<\varepsilon^2\}}\cap {\{|y-t|<\varepsilon^2\}}\cap {\{|x-s|<\varepsilon^2\}}\cap{\{|x-u|<\varepsilon^2\}}$$
and
$$B:={\{|x-s|<\varepsilon^2\}}\cap{\{|y-t|<\varepsilon^2\}}\cap{\{|z-t|<\varepsilon^2\}}\cap{\{|w-v|<\varepsilon^2\}}.$$

Integrating with respect to $u,s,t,v$ in  the term corresponding to
$I_A$  and with respect to $x,y,z,w$ in the term corresponding to
$I_B$, we obtain the desired result.

The proof of Theorem \ref{nouteo} is now complete. \hfill $\Box$

\end{dem}

\begin{enumerate}

\item[]{F. Avram,}{ Weak convergence of the variations, iterated integrals and Dol\'{e}ans-Dade exponentials of sequences of
semimartingales. }{Ann. Probab., 16(1) (1988) 246-250}

\item[]{ X. Bardina and M. Jolis,} { Weak convergence
to the multiple Stratonovich integrals.} {Stochastic Process. Appl.,
90(2) (2000) 277-300}

\item[]{P. Billinsgley,}{ Convergence of probability
measures.} {John Wiley and Sons} (1968)

\item[]{K. It\^{o},}{ Multiple Wiener Integral.} {J. Math. Soc. Japan, 3 (1951), 157-169}

\item[]{M. Kac,}{ A stochastic model related to the telegrapher's equation. Reprinting of an article published in
1956.} {Rocky Mountain J. Math., 4 (1974), 497-509}


\item[]{J.L. Sol\'e and F. Utzet,}{ Stratonovich
integral and trace.} {Stochastics Stochastics Rep., 29(2) (1990)
203-220 }

\item[]{ D. Stroock,} {Topics in Stochastic
Differential Equations.} {(Tata Institute of Fundamental Research,
Bombay).} {Springer Verlag} (1982)
\end{enumerate}
\end{document}